\documentclass[reqno,11pt]{amsart}

\usepackage{amsfonts,amsmath,amssymb,fullpage,graphicx}
\usepackage{verbatim}

\newcommand{\de}{\mathrm{d}}
\newcommand{\e}{\mathrm{e}}
\newcommand{\nn}{\nonumber}
\newcommand{\ug}{\!\!\!\!&=&\!\!\!\!}

\begin{document}

\makeatletter
\title{Integrals of Bessel functions}
\author{D. Babusci}
\address{INFN - Laboratori Nazionali di Frascati, via  E. Fermi, 40, IT 00044 Frascati (Roma), Italy}
\email{danilo.babusci@lnf.infn.it}
\author{G. Dattoli}
\address{ENEA - Centro Ricerche Frascati, via  E. Fermi, 45, IT 00044 Frascati (Roma), Italy}
\email{giuseppe.dattoli@enea.it}
\author{B. Germano}
\address{Dipartimento di Metodi e Modelli Matematici per le Scienze Applicate, 
Universit\`a degli Studi di Roma ``Sapienza", via A. Scarpa, 14, Roma, Italy}
\email{germano@dmmm.uniroma1.it}
\author{M. R. Martinelli}
\address{Dipartimento di Metodi e Modelli Matematici per le Scienze Applicate, 
Universit\`a degli Studi di Roma ``Sapienza", via A. Scarpa, 14, Roma, Italy}
\email{martinelli@dmmm.uniroma1.it}
\author{P. E. Ricci}
\address{International Telematic University UniNettuno, 
Corso Vittorio Emanuele II, 39, 00186 Roma, Italy}
\email{paoloemilioricci@gmail.com}


\subjclass[2000]{Primary 33}

\begin{abstract}
We use the operator method to evaluate a class of integrals involving Bessel or Bessel-type functions. 
The technique we propose is based on the formal reduction of these family of functions to Gaussians.
\end{abstract}

\maketitle

Integrals of the type
\begin{equation}
\label{eq:Inab}
I_n (a, b, \alpha) = \int_{- \infty}^\infty \de x \,(a\,x + b)^n\,\e^{- \alpha\,x^2}\,
\end{equation}
can be calculated using a general procedure based on the method of generating function (GF) \cite{Babusci}. Multiplying both sides of 
eq. \eqref{eq:Inab} by $t^n /n!$ and summing up over $n$, we find 
\begin{eqnarray}
\label{eq:Gab}
G (a, b, \alpha) \ug \sum_{n = 0}^\infty \frac{t^n}{n!}\,I_n (\alpha, \beta) = \e^{\,b\,t}\,\int_{- \infty}^\infty \de x \,\e^{- \alpha\,x^2 + t\,a\,x} \nn \\
\ug \sqrt{\frac{\pi}{\alpha}}\,\exp \left(\frac{t^2\,a^2}{4\,\alpha} + b\,t\right)\,.
\end{eqnarray}
By exploiting the GF of the two-variable Hermite polynomials \cite{Appell}
\begin{equation}
\label{eq:Herm}
\sum_{n = 0}^\infty \frac{t^n}{n!}\,H_n (x,y) = \e^{x\,t + y\,t^2} \qquad\qquad H_n (x,y) = n!\,\sum_{k = 0}^{[n/2]} \frac{x^{n - 2\,k}\,y^k}{(n - 2\,k)!\,k!}\,,
\end{equation}
from eq. \eqref{eq:Gab} we get 
\begin{equation}
G (a, b, \alpha) = \sqrt{\frac{\pi}{\alpha}}\,\sum_{n = 0}^\infty \frac{t^n}{n!}\,H_n \left(b, \frac{a^2}{4\,\alpha}\right)\,, \nn
\end{equation}
and, thus 
\begin{equation}
\label{eq:InHn}
I_n (a, b, \alpha) = \sqrt{\frac{\pi}{\alpha}}\,H_n \left(b, \frac{a^2}{4\,\alpha}\right)\,.
\end{equation}

In a recent series of papers \cite{Babusci} it has been proved that the Bessel functions can be formally reduced to Gaussians according to the 
following representation
\begin{equation}
J_\nu (x) = \left(\frac{x\,\hat{c}}2\right)^\nu\, \exp\left\{- \hat{c}\,\left(\frac{x}2\right)^2\right\}\,\varphi (0)
\end{equation}
where
\begin{equation}
\label{eq:cop}
\hat{c}^{\,\mu}\,\varphi (0) = \varphi (\mu)\,, \qquad\qquad \varphi (\mu) = \frac1{\Gamma (\mu + 1)}\,.
\end{equation}
The properties of the operator $\hat{c}$ have been extensively discussed in \cite{Babusci}. By treating it as an ordinary constant, we find, 
for example, 
\begin{equation}
B_0 (\alpha) = \int_{- \infty}^\infty \de x\, J_0 (\sqrt{\alpha}\,x) = \frac2{\sqrt{\alpha}}\,, 
\end{equation}
and, in analogy with eq. \eqref{eq:Inab}, for $\nu > n$ ($\nu \in \mathbb{R}$)  
\begin{eqnarray}
\int_{- \infty}^\infty \de x\, J_\nu (\sqrt{\alpha}\,x)\,\frac{(a\,x + b)^n}{(\sqrt{\alpha}\,x)^\nu} 
\ug \frac1{2^{\nu - 1}}\,\frac{\sqrt{\pi}}{\sqrt{\alpha}}\,\hat{c}^{\,\nu - 1/2}\,H_n \left(b,\frac{a^2}{\alpha\,\hat{c}}\right) \nn \\
\ug \frac1{2^{\nu - 1}}\,\frac{\sqrt{\pi}}{\sqrt{\alpha}}\,B_n \left(b, \frac{a^2}{\alpha}; \nu\right) 
\end{eqnarray}
where we have introduced the polynomials
\begin{equation}
\label{eq:Bn}
B_n (x, y; \nu) =  n!\,\sum_{k = 0}^{[n/2]} \frac{x^{n - 2\,k}\,y^k}{(n - 2\,k)!\,k!\,\Gamma (\nu - k + \frac12)}\,,
\end{equation}
whose properties will be briefly described later. As a consequence of this result, for any function $f (x)$ that can be written as 
\begin{equation}
f (x, m) = \sum_{k = 0}^{m} f_k\,x^{k} \qquad\qquad (m < \nu)\,,
\end{equation}
the following identify holds  
\begin{equation}
\int_{- \infty}^\infty \de x\, J_\nu (\sqrt{\alpha}\,x)\, \frac{f (a\,x + b)}{(\sqrt{\alpha}\,x)^\nu} = 
\frac1{2^{\nu - 1}}\,\frac{\sqrt{\pi}}{\sqrt{\alpha}}\,\sum_{k = 0}^{m} f_k\,B_k \left(b, \frac{a^2}{\alpha}; \nu\right)\,.
\end{equation}

It is interesting to note that the linear combination of Bessel functions
\begin{equation}
F_n (x; a, b) = \sum_{k = 0}^n \binom{n}{k}\,a^{\,n - k}\,b^k\,J_{n - k} (x) 
\end{equation}
can be written as
\begin{equation}
F_n (x; a, b) = \left(\frac{a\,x}2\,\hat{c} + b\right)^n\,\exp\left\{- \hat{c}\,\left(\frac{x}2\right)^2\right\}\,\varphi (0)
\end{equation}
and, therefore,  according to eq. \eqref{eq:InHn}, the following identity holds 
\begin{eqnarray}
\int_{- \infty}^\infty \de x\,F_n (x; a, b) \ug 2\,\sqrt{\pi}\,\hat{c}^{- 1/2}\,H_n \left(b, \frac{a^2}4\,\hat{c}\right)\,\varphi (0) \nn \\
\ug 2\,\sqrt{\pi}\,n!\,\sum_{k = 0}^{[n/2]} \frac{b^{\,n - 2\,k}\,a^{\,2\,k}}{4^k\,(n - 2\,k)!\,k!\,\Gamma \left(k + \frac12\right)}\,.
\end{eqnarray}

The method illustrated is more general than it may appear and can indeed be extended to other families of Bessel-like functions.  
For example, the spherical Bessel functions \cite{Andrews} 
\begin{equation}
j_n (x) = \sqrt{\frac{\pi}{2\,x}}\,J_{n + 1/2} (x)
\end{equation}
within the present formalism can be written as
\begin{equation}
\label{eq:jn}
j_n (x) = \frac{\sqrt{\pi}}{2^{\,n +1}}\,x^n\,\hat{c}^{\,n + 1/2}\,\exp\left\{- \hat{c}\,\left(\frac{x}2\right)^2\right\}\,\varphi (0)\,.
\end{equation}
By using this expression, the integral 
\begin{equation}
b_n = \int_{- \infty}^\infty \de x\,j_n (x)\,,
\end{equation}
can easily be calculated by exploiting the GF method, getting
\begin{equation}
b_{2\,n} = \frac{\pi}{2^{2\,n}}\,\frac{(2\,n)!}{(n!)^2}\,, \qquad\qquad b_{2\,n + 1} = 0\,.
\end{equation}

The Struve functions are defined by the series \cite{Andrews}
\begin{equation}
\mathbf{H}_\nu (x) = \sum_{k = 0}^\infty \frac{(- 1)^k}{\Gamma (k + \frac32)\,\Gamma (k + \nu + \frac32)}\,
\left(\frac{x}2\right)^{2\,k + \nu + 1}
\end{equation}
Albeit in a slightly different form, even in this case we can apply the operational method, and this functions can be 
written as
\begin{equation}
\mathbf{H}_\nu (x) = \hat{c}_1^{\,1/2}\,\hat{c}_2^{\,\nu + 1/2}\,\left(\frac{x}2\right)^{\,\nu + 1}\,
\frac1{1 + \hat{c}_1\,\hat{c}_2\,\left(\frac{x}2\right)^2}\,\varphi_1 (0)\,\varphi_2 (0)
\end{equation}
where the operator $\hat{c}_i$, ($i = 1,2$) acts only on $\varphi_i$ and verifies the identity \eqref{eq:cop}. Moreover, 
the use of the Laplace allows us to write
\begin{equation}
\mathbf{H}_\nu (x) =  \hat{c}_1^{\,1/2}\,\hat{c}_2^{\,\nu + 1/2}\,\left(\frac{x}2\right)^{\,\nu + 1}\,
\int_0^\infty \de s\,\exp\left\{- s\,\left[1 + \hat{c}_1\,\hat{c}_2\,\left(\frac{x}2\right)^2\right]\right\}\,\varphi_1 (0)\,\varphi_2 (0)
\end{equation}
that can be used to prove that, for $\mu + \nu$ not an even integer, the following identity holds 
\begin{equation}
\int_0^\infty \de x\,x^\mu\,\mathbf{H}_\nu (x) = (- 1)^{\,\mu + \nu}\,\frac{2^{\,\mu}\,\pi}{\sin \left[(\mu + \nu)\,\frac{\pi}2\right]}\,
\frac1{\Gamma \left(\frac{1 - \mu - \nu}2\right)\,\Gamma \left(\frac{1 - \mu + \nu}2\right)}\,.
\end{equation}

As a final example, we consider the Wright-Bessel functions \cite{Sriva}
\begin{equation}
W_{\alpha, \beta} (x) = \sum_{k = 0}^\infty \frac{x^k}{k!\,\Gamma (k\,\alpha + \beta)} 
\end{equation}
($\beta \notin \mathbb{Z}$ when $\alpha < 0$) that can formally be defined as 
\begin{equation}
\label{eq:WB}
W_{\alpha, \beta} (x) = \hat{c}^{\,\beta - 1}\,\exp\left(\hat{c}^\alpha \,x\right)\,\varphi (0)\,.
\end{equation}
By using this expression it's easy to show that 
\begin{equation}
\int_{- \infty}^\infty \de x\,W_{\alpha, \beta} (- x^2) = \sqrt{\pi}\,\frac1{\Gamma (\beta - \alpha/2)}\,,
\end{equation}
and ($d > 0$) 
\begin{equation}
\int_0^\infty \de x\,W_{\alpha, \beta} (- x)\,\e^{- d\,x} = \frac1d\,E_{\alpha, \beta} \left(- \frac{1}{d}\right)
\end{equation}
where 
\begin{equation}
E_{\alpha, \beta} (x) = \sum_{k = 0}^\infty \frac{x^k}{\Gamma (\alpha\,k + \beta)}
\end{equation}
is the modified Mittag-Leffler function \cite{Sriva}.

Now, let us coming back to the polynomials $B_n (x, y; \nu)$, introduced in eq. \eqref{eq:Bn}. These are auxiliary Hermite-like 
polynomials, and their properties can be deduced from those of Hermite polynomials written in umbral form. For example, 
by using the first of eq. \eqref{eq:Herm} and eq. \eqref{eq:WB}, it can be shown that the generating function is given by  
\begin{equation}
\sum_{n = 0}^\infty \frac{t^n}{n!}\,B_{n} (x, y; \nu) = \e^{x\,t}\,W_{-1, \nu + 1/2} (y\,t^2) 
\end{equation}
It is also interesting to note that these polynomials satisfy the differential equation
\begin{equation}
\hat{D}_y\,B_{n} (x, y; \nu) = \partial_x^2\,B_{n} (x, y; \nu) \qquad\qquad B_{n} (x, 0; \nu) = x^n\,\varphi (\nu - \frac12)
\end{equation}
where the derivative operator $\hat{D}_y = \hat{c}\,\partial_y$ has been introduced. Therefore, we get 
\begin{equation}
B_{n} (x, y; \nu) = \exp\left(\hat{c}^{- 1}\,y\,\partial_x^2\right)\,x^n\,\varphi (\nu - \frac12)\,,
\end{equation}
that, taking into account the operational definition of the Hermite polynomials
\begin{equation}
H_n (x, y) = \exp\left(y\,\partial_x^2\right)\,x^n\,,
\end{equation}
can be easily shown to coincide with the expression given in eq. \eqref{eq:Bn}. 
These polynomials can be framed within the Appell family. Nevertheless, 
they deserve to be studied carefully and this will be done in a next paper.

\end{document}